\documentclass[]{interact}

\usepackage{epstopdf}% To incorporate .eps illustrations using PDFLaTeX, etc.
\usepackage[caption=false]{subfig}% Support for small, `sub' figures and tables

\usepackage[numbers,sort&compress]{natbib}% Citation support using natbib.sty
\bibpunct[, ]{[}{]}{,}{n}{,}{,}% Citation support using natbib.sty
\makeatletter% @ becomes a letter
\def\NAT@def@citea{\def\@citea{\NAT@separator}}% Suppress spaces between citations using natbib.sty
\makeatother% @ becomes a symbol again

\theoremstyle{plain}% Theorem-like structures provided by amsthm.sty
\newtheorem{theorem}{Theorem}[section]
\newtheorem{lemma}[theorem]{Lemma}

\theoremstyle{definition}

\theoremstyle{remark}
\newtheorem{remark}{Remark}

\begin{document}

\articletype{ARTICLE TEMPLATE}% Specify the article type or omit as appropriate

\title{ Local error estimates and post processing for the Galerkin boundary element method on polygons}

 \author{  Dedicated to Prof. George C. Hsiao on the occasion of his 90th birthday.  \\[0.5cm]
\name{ Thomas Hartmann  \textsuperscript{a} \thanks{CONTACT
 Thomas Hartmann: Email: thomas.hartmann@thu.de} 
and   Ernst~P. Stephan \textsuperscript{b}}
\affil{  \textsuperscript{a}  University of Applied Sciences Ulm, Germany;  \textsuperscript{b} Institut  f\"ur Angewandte Mathematik,  Leibniz University  Hannover, Germany }}

\maketitle

\begin{abstract}
In this paper we give local error estimates in Sobolev norms
for the Galerkin method applied to strongly elliptic pseudodifferential equations
on a polygon.   By using the K-operator, an operator which averages
the values of the Galerkin solution, we construct improved approximations. 
\end{abstract}


\section{Introduction}  \label{lfa}
In the following we consider strongly elliptic integral equations on a polygon.
It is well known that the behaviour of the 
approximation depends heavily on the characteristics of the settings of 
the problem. Non-smoothness of domains and singularities of given data 
produce singularities of the exact solution, and reduce the global rate of 
convergence of the approximations. It is then natural to ask whether the 
accuracy of the approximation is better in regions where the exact solution
is smooth.

For strongly elliptic integral equations on a smooth and closed curve in 
$\mathbb{R}^2$, local error estimates were proved by Saranen \cite{sar}. His proof
was then modified by Tran to derive local estimates on open arcs \cite{thanh}.
The results were extended to open surfaces in $\mathbb{R}^3$ \cite{thste}. 
In the present paper for polygons local estimates are shown together
with corresponding numerical results. 
Recently in \cite{faust} local estimates of the bem were obtained also for polyhedral domains.\\
We consider in Section \ref{lrg} first the the local regularity 
of the Galerkin solution and then we derive local error estimates. 
In Section \ref{nummm} numerical results are presented, which confirm the 
theoretical results. The K-operator is investigated in Section 
\ref{kop}. 

Let $\Gamma$ be a polygon. For a pseudodifferential operator $A$ of 
order $2\alpha$ \cite{see}\cite{gwinnerstephan} 
we consider  the equation
\begin{equation} \label{ia}
A\psi=f
\end{equation}
assuming that $A:H^\alpha(\Gamma)\to H^{-\alpha}(\Gamma)$ is a continuous,
bijective mapping acting between the Sobolev spaces $H^s(\Gamma), s\in \mathbb{R}$.
These function spaces are introduced for an arbitrary subsegment $\Gamma'$ 
of $\Gamma$ 
$$H^s(\Gamma'):=\{ u\vert_{\Gamma'} : u \in H^{s}(\Gamma)\} \qquad s\in \mathbb{R}.
$$
Furthermore we assume that
there exists a constant  $\gamma>0$, such 
that for all $\psi \in H^\alpha (\Gamma)$ there holds
\begin{equation} \label{sel}
 \langle \;A\psi , \psi \;\rangle \geq \gamma \Vert \psi \Vert_\alpha^2.
\end{equation}
where $ \langle , \rangle $ denotes the duality between $ H^{-\alpha}(\Gamma)$ 
and $ H^{\alpha}(\Gamma)$. 
On a quasiuniform mesh $\triangle=\{x_1,...,x_n\}$,
$x_i \in \Gamma, i=1,...,n$ we define smoothest splines of order $r$ 
(or degree $r-1$)dot with breakpoints in $\triangle$. The space of these splines is 
denoted by $S^r_h(\Gamma)$, where $h$ is the maximum value of the step sizes. 
This space will be used both as test and trial space in the Galerkin 
scheme.

We  
assume the stability condition: there exists a constant $c>0$ such that 
for $0<h<h_0$
\begin{equation} \label{bb}
\inf_{\phi \in S_h^r} \sup_{\psi \in S_h^r} 
\frac{\vert \langle A\phi,\psi\rangle \vert}
{\Vert \phi \Vert_{\alpha} \Vert \psi \Vert_{\alpha} }
\geq  c >0.
\end{equation}

Now for $\Gamma_0 \subset \Gamma$ we are concerned with the local Galerkin error 
$\Vert \psi - \psi_h \Vert_{H^t(\Gamma_0)}$ between the solution 
$\psi$ of (\ref{ia}) and its Galerkin solution $\psi_h$  
which satisfies
\\
{\it Find $ \psi_h \in S^r_h(\Gamma)\subset  H^\alpha(\Gamma)$ such that for all 
$\phi \in S^r_h(\Gamma)$ 
\begin{equation} \label{galeq}
 \langle \;A\psi_h , \phi \;\rangle = \langle \;f,\phi \;\rangle .
\end{equation} }

 \section{Local error estimates} \label{lrg}
 Let $\Gamma^j$ denote a fixed straight line segment of the polygon $\Gamma=\cup_{j=1}^J \Gamma^j$ and let $\beta_j$ denote the interior
angles between  $\Gamma^j$ and $\Gamma^{j+1}$.  Further we use nested 
subsegments  $\Gamma_i,i=0,...,M,$ $\Gamma_*$ and $\Gamma^*$ of $\Gamma^j$ 
and associated 
cut-off functions $\omega_i$
$$
\Gamma_* \subset\subset\Gamma_0 \subset\subset\Gamma_1 \subset\subset \ldots  
\subset\subset\Gamma_M \subset\subset \Gamma^* \subset\subset \Gamma^j
$$
\begin{equation} \label{cut}
\omega_i = 1 \mbox{ on } \Gamma_i \qquad
\omega_i \in C^\infty(\Gamma) \qquad
supp \;\omega_i \subset\subset \Gamma_{i+1}.
\end{equation}
Here $X \subset\subset Y$ means the closure of $X$  is contained in the interior of $Y$. 

For abbreviation we use the notations
$$\qquad \Vert \cdot \Vert_{s} := \Vert \cdot \Vert_{H^s(\Gamma)}
 \qquad \Vert \cdot \Vert_{s,\Gamma_j} := \Vert \cdot \Vert_{H^s(\Gamma_j)}
\qquad \Vert \cdot \Vert_{\Gamma_j} = \Vert \cdot \Vert_{0,\Gamma_j}.$$

 Let $\Gamma^j$ be an arbitrary but 
fixed segment and $\tilde{\Gamma}$ a smooth closed curve containing $\Gamma^j$.
In the following we derive local error estimates for a subsegment 
$\Gamma_0$, which is contained in the straight line segment $\Gamma^j$. 
Furthermore let the pseudodifferential operator  $\tilde{A}$ be defined 
as $A$ in (\ref{ia})
but on $\tilde{\Gamma}$ instead of on ${\Gamma}$.
Now we recall some basic properties of pseudo-differential operators:

\begin{lemma} \label{com} \cite{see}
Let $\omega $ be a $C^\infty$ cut-off function. Then for any $s \in \mathbb{R}$ 
we have for  $\phi\in H^{s-2}(\tilde{\Gamma})$
\begin{equation} \label{comest}
\Vert (\omega \tilde{A} - \tilde{A} \omega) \phi \Vert_{H^s(\tilde{\Gamma})}
\leq C\,  \Vert \phi \Vert_{H^{s-2}(\tilde{\Gamma})}
\end{equation}	
\end{lemma}

\begin{lemma} \label{gla} \cite{see}
Let $\omega$ and $\tilde{\omega}$ be cut-off functions with disjoint supports.
Then we have for arbitrary but fixed $\beta>\alpha$ and for all $\phi\in H^{\beta}(\tilde{\Gamma})$.
\begin{equation} \label{glaest} 
\Vert \tilde{\omega} \tilde{A} \omega \phi \Vert_{H^{-\alpha}(\tilde{\Gamma})}
\leq C\,  \Vert \phi \Vert_{H^{\beta}(\tilde{\Gamma})}
\end{equation}
\end{lemma}

Higher convergence rates of the Galerkin scheme are expected if the exact solution of the integral
equation (\ref{ia}) is smoother in some subsegments. Therefore 
it is worthwhile to study the local regularity of the solution.

\begin{theorem} \label{ntn}
Let  $f\in H^{s-2\alpha}(\Gamma^*)\cap H^{-\alpha}(\Gamma)$ with  
$s\geq \alpha$. Then there holds for the solution 
$\psi \in H^\alpha(\Gamma)$ of (\ref{ia})
\begin{equation} \label{reg}
\Vert \psi \Vert_{H^{s}(\Gamma_*)} \leq C\, \{ \Vert f \Vert_{H^{s-2\alpha}(\Gamma^*)}
+ \Vert \psi \Vert_{H^{\alpha}(\Gamma)} \}
\end{equation}
\end{theorem}
Proof: With cut-off functions $\omega_i$ as in (\ref{cut}) we have
$\omega_0 A \psi = \omega_0 f$
and further
\begin{equation} \label{dh}
\omega_0 A \omega_1 \psi
=- \omega_0 A (1-\omega_1) \psi+  \omega_0 f
=: \omega_0 h.
\end{equation}
Since $\omega_0$ and $ (1-\omega_1)$ have disjoint support we have for any 
$s> \alpha$
$$
\Vert \omega_0 A (1-\omega_1) \psi \Vert_{H^{s-2\alpha}(\Gamma)}
\leq C\,  \Vert \psi \Vert_{H^{\alpha}(\Gamma)}.
$$
This yields with  (\ref{dh})
$$ \Vert \omega_0 h \Vert_{H^{s-2\alpha}(\Gamma)} \leq
\Vert \omega_0 A (1-\omega_1) \psi \Vert_{H^{s-2\alpha}(\Gamma)}
+ \Vert \omega_0 f \Vert_{H^{s-2\alpha}(\Gamma)}
$$ 
\begin{equation} \label{nh}
\leq C\, \{\Vert \psi \Vert_{H^{\alpha}(\Gamma)}+\Vert f\Vert_{H^{s-2\alpha}(\Gamma_1)}\}.
\end{equation}
Since $supp\; \omega_0$ and $supp\; \omega_1$ are contained in $ \Gamma^j$ 
we have:
$$ \omega_0 A          \omega_1 \psi =
   \omega_0 \tilde{A}  \omega_1 \psi
$$
and from (\ref{dh}) we obtain after multiplying by   $\tilde{A}^{-1}$
$$ \tilde{A}^{-1} \omega_0 h
=  \tilde{A}^{-1} \omega_0 \tilde{A} \omega_1 \psi
= \tilde{A}^{-1} \omega_0 \tilde{A} \omega_1 \psi
- \omega_0 \tilde{A}^{-1} \tilde{A} \omega_1 \psi
+ \omega_0  \omega_1 \psi
$$
\begin{equation} \label{th}
= \left[\tilde{A}^{-1} \omega_0- \omega_0 \tilde{A}^{-1}\right]
\tilde{A} \omega_1 \psi
+ \omega_0 \psi.
\end{equation}
Now  $\tilde{A}^{-1}$  is a   pseudo-differential operator of order 
$-2\alpha$ und hence the commutator 
$\left[\tilde{A}^{-1} \omega_0- \omega_0 \tilde{A}^{-1}\right]$
is a pseudo-differential operator of order $-2\alpha-1$. From (\ref{th}) 
we have
\begin{eqnarray*}
\Vert \omega_0 \psi \Vert_{H^{s}(\Gamma)} 
& \leq &
 \Vert \tilde{A}^{-1}\omega_0 h \Vert_{H^{s}(\Gamma)}
+\left\Vert \left[\tilde{A}^{-1} \omega_0-
             \omega_0 \tilde{A}^{-1}\right]
\tilde{A} \omega_1 \psi\right\Vert_ {H^{s}(\Gamma)}
\\
& \leq C&  \{\Vert \omega_0 h \Vert_{H^{s-2\alpha}(\Gamma)} +
     \Vert \omega_1 \psi \Vert_{H^{s-1}(\Gamma)} \}
\end{eqnarray*}
and with (\ref{nh})
$$
\Vert \omega_0 \psi \Vert_{H^{s}(\Gamma)} \leq
C\, \{ \Vert \psi \Vert_{H^{\alpha}(\Gamma)} +
    \Vert f  \Vert_{H^{s-2\alpha}(\Gamma_1)}+
     \Vert \omega_1 \psi \Vert_{H^{s-1}(\Gamma)} \}.
$$
Repeating the argument we have
$$
\Vert \omega_1 \psi \Vert_{H^{s-1}(\Gamma)} \leq
C\, \{ \Vert \psi \Vert_{H^{\alpha}(\Gamma)} +
    \Vert f  \Vert_{H^{s-2\alpha-1}(\Gamma_2)}+
     \Vert \omega_2 \psi \Vert_{H^{s-2}(\Gamma)} \}
$$
and taking more subsegments yields the desired estimate.\qed

Now we change the notation and use nested subsegments  
$\Gamma_i,i=0,...,M,$ 
and associated cut-off functions $\omega_i$
$$
\Gamma_0 \subset\subset\Gamma_1 \subset\subset \ldots  
\subset\subset\Gamma_M \subset\subset \Gamma_* 
\subset\subset \Gamma^* \subset\subset \Gamma^j$$
defined as in (\ref{cut}).
In order to analyse the local error of the Galerkin scheme
we use the following spline spaces for $j=0,1,...,M$
$$
{S_h^{r,0}}(\Gamma_j)=\{ \phi \in S^r_h : \mbox{supp} \phi \subset \Gamma_j\}
$$ $$
S^r_h(\Gamma_j)=\{ v\in H^{-1/2}(\Gamma) : 
                v\vert_{\Gamma_j} = \phi \vert_{\Gamma_j} \mbox{ for some }
                \phi \in S_h^r\}
$$

Due to Theorem \ref{ntn} there holds
\begin{equation} \label{a1}
\psi \in H^{\alpha}(\Gamma) \cap H^s(\Gamma_*) \qquad 
 \mbox{for }\alpha \leq s 
\end{equation}
\begin{equation} \label{a2}
\psi_h \in S^r_h(\Gamma_*) \qquad \mbox{with} \quad r\geq s
\end{equation}
\begin{equation} \label{a3}
\langle \; A(\psi -\psi_h),\phi \; \rangle =0\qquad
\mbox{for all } \phi \in S^{r,0}_h(\Gamma_*).
\end{equation}
Note that in (\ref{a2}) we only assume  $\psi_h$ is a spline 
on the subsegment $\Gamma_*$ and can be extended arbitrarily onto $\Gamma$
provided the extension belongs to $H^{\alpha}(\Gamma)$. 
Saranen assumed in \cite{sar} that $\psi_h$ in (\ref{a2}) 
is a spline on the whole boundary $\Gamma$. 
\\
In the following we show an estimate for the local Galerkin error for 
the solution of the pseudodifferential equation (\ref{ia}) on a polygon.

\begin{theorem} \label{theo}
Let $\psi \in H^s(\Gamma_*)$ be the solution of (\ref{ia}) and $\psi_h$ satisfy
(\ref{a2}) and (\ref{a3}). Then for any $\beta <\alpha$ 
we have for $-r +2\alpha \leq t\leq s \leq r$ and $t<r+\alpha$: there exists
a constant $C=C(\beta) $ independent of $h, \psi$ and $\psi_h$
\begin{equation} \label{inn}
 \Vert \psi - \psi_h \Vert_{t,\Gamma_0}
\leq C\,  \{h^{s-t} \Vert \psi \Vert_{s,\Gamma_*}+
h^\sigma\Vert \psi- \psi_h \Vert_{\beta,\Gamma} \} 
\end{equation}
where $\sigma =0$  for $t\leq \alpha$ and $\sigma=\alpha-t$  for $\alpha< t$.
\end{theorem}

In order to derive the local error estimate of Theorem \ref{theo} we 
analyse (\ref{a3}) and use the following technical results 
(Lemma \ref{app} -\ref{sup}).
The spline spaces defined above have the following properties \cite{bab,sar}.

\begin{lemma} \label{app} Assume that $t_0< r-1/2$ and $q\in \mathbb{N}$. 
Let $\phi \in H^s(\Gamma)$ with $t_0\leq s \leq r$.\\
(a) There exists $\xi \in S^r_h$ such that for all $t<t_0$:
$$\Vert  \phi - \xi \Vert_t 
\leq C h^{s-t} \Vert \phi \Vert_{s,\Gamma}
$$
(b) For any $j=0,..., M-2$, there exists $h_0>0$ and $\xi \in {S^{r,0}_h(\Gamma_{j+2})}$ such that
$$\Vert \omega_j \phi - \xi \Vert_t 
\leq C h^{s-t} \Vert \phi \Vert_{s,\Gamma_{j+1}},
$$
for all $t \in [-q,t_0]$ and $0<h\leq h_0$.
\end{lemma}

\begin{lemma} \label{inv} 
For all $t\leq s <r-1/2$, $j=0,...,M-1$, and $\phi \in S^r_h$, there holds
\begin{eqnarray}
 \Vert \phi \Vert_s & \leq&  C h^{t-s} \Vert \phi \Vert_t \nonumber \\
 \label{locinv}
 \Vert \omega_j \phi \Vert_s 
& \leq & C h^{t-s} \Vert \phi \Vert_{t,\Gamma_{j+1}} 
\end{eqnarray}
\end{lemma}

\begin{lemma} \label{sup} Assume $t_0<r-1/2$ and $q\in \mathbb{N} $. Then for any $\phi
\in S_h^r$ and $j=0,...,M-2$ there exists $h_0>0$ and $\xi \in 
{S_{h}^{r,0}}(\Gamma_{j+2})$ such that
$$\Vert \omega_j \phi - \xi \Vert_t 
\leq C h^{s+1-t} \Vert \phi \Vert_{s,\Gamma_{j+1}},
$$
for all $t \in [-q,t_0], s\leq r-1$ and $0<h\leq h_0$.
\end{lemma}
Now we prove the local error estimate given in Theorem \ref{theo}.\\
{\bf Proof ot Theorem \ref{theo}:}
We split the local error into three terms (see (\ref{decom}))
which are estimated seperately in Lemma \ref{dec1} - \ref{dec3}. Combining the results of these lemmas we prove with Lemma \ref{erinen} the estimate (\ref{inn}) for $ t=\alpha $. Then the proof 
of Theorem \ref{theo} can be completed by following the arguments in \cite{thste}. For brevity we
omit the corresponding details. \qed

Let us rewrite (\ref{galeq}) as

{\it For any $\psi \in H^{\alpha}(\Gamma)$, find 
$ G\psi \in S^r_h$ such that 
\begin{equation} \label{galaux}
 \langle \;A(\psi -G\psi) , \phi \;\rangle = 0
\mbox{  for any } \phi \in S^r_h
\end{equation} }

To simplify the notation we introduce 
 $e:=\psi-\psi_h$ and $\tilde{\phi}:=\omega_0 \phi$
for arbitrary $\phi \in H^s(\Gamma)$. We decompose the 
local error $\tilde{e}:=\omega_0 e$ 
in
\begin{equation} \label{decom}
\tilde{e} =(\tilde{\psi} -  G\tilde{\psi}) 
                         + (G\tilde{\psi} - G\tilde{\psi_h})
                                          +(G\tilde{\psi_h} - \tilde{\psi_h}),
\end{equation}
and estimate each of the terms in parantheses seperately.

\begin{lemma} \label{dec1}
Under the assumption (\ref{a1}) for $\psi$ and $\tilde{\psi}:=\omega_0 \psi$ 
we have for $0<h<h_0$ and $\alpha<s<r$:
$$\Vert \tilde{\psi} -  G\tilde{\psi} \Vert_{\alpha} 
\leq C\,  h^{s-\alpha} \Vert \psi \Vert_{s,\Gamma_*}$$
\end{lemma}
Proof: Due to (\ref{bb}) the Galerkin method yields quasioptimal error
estimates.
Using  Lemma \ref{app} we have the estimate:
$$ \Vert \tilde{\psi} -  G\tilde{\psi} \Vert_{\alpha} 
\leq C\,  \inf_{\phi \in S_h^r} \Vert \phi - \tilde{\psi}\Vert_{\alpha}
\leq C\,  h^{s-\alpha} \Vert \tilde{\psi} \Vert_{s} 
\leq C\,  h^{s-\alpha} \Vert \psi         \Vert_{s,\Gamma_*}. $$
This gives the proof. \mbox{\qed}
\vspace{0.2cm}

\begin{lemma} \label{dec2}
 Under the assumptions (\ref{a1}-\ref{a3}) for any fixed $\beta$ with 
$\beta <\alpha$ there holds
$$\Vert G\tilde{\psi} - G\tilde{\psi_h} \Vert_{\alpha} \leq C\,  \left\{
h \Vert e \Vert_{\alpha,\Gamma_*} 
+ \Vert e \Vert_{\alpha-1,\Gamma_*}
+ \Vert e \Vert_{\beta} \right\}$$
\end{lemma}
Proof: With (\ref{bb}) we obtain 
$$\Vert G\tilde{e}\Vert_{\alpha}\leq C\,  \sup_{\phi \in S^r_h} 
\frac{ \vert \langle A G\tilde{e}, \phi \rangle \vert}{
\Vert\phi\Vert_{\alpha}}$$
and further with (\ref{galaux}) and $\omega_M$ from (\ref{cut}):
$$ \langle \; A G \tilde{e} , \phi \;\rangle =
\langle \; A  \tilde{e} , \phi \;\rangle = 
\langle \;\omega_M A\tilde{e},\phi\; \rangle +
\langle \;(1-\omega_M) A\tilde{e},\phi\; \rangle.
$$
We estimate the second term with Lemma \ref{gla} and decompose the first term
$$
| \langle \;(1-\omega_M) A\tilde{e},\phi\; \rangle |
\leq \Vert (1-\omega_M) A\tilde{e}\Vert_{-\alpha}
\Vert\phi\Vert_{\alpha} 
\leq C\,  \Vert e \Vert_{\beta} \Vert\phi\Vert_{\alpha}
$$
$$
 \langle \; A \tilde{e}, \omega_M \phi \; \rangle =
\langle \; \omega_0 A e, \omega_M \phi \; \rangle +
\langle \; A \tilde{e}-\omega_0 A e, \omega_M \phi \; \rangle =:B+D
$$
Let $\xi \in {S_{h}^{r,0}}(\Gamma_2)$ such that (see Lemma \ref{sup} )
\begin{equation} \label{supbsp} 
\Vert\omega_0\phi - \xi \Vert_{\alpha} 
\leq C\,   h \Vert\phi\Vert_{\alpha,\Gamma_{j+1}}
\leq C\,   h \Vert\phi\Vert_{\alpha},
\end{equation}
here we used $r\geq \alpha+1$. If this assumption is not fulfilled 
the estimate is still valid with exchanging $h$ by $h^\tau$, where
$\tau:= \min \{ 1, r-\alpha \}$.
With (\ref{a3}) and supp$ \;(\omega_0\phi - \xi) \subset\subset \Gamma_2$ 
we have
\begin{eqnarray*}
 B& =& \langle \;Ae , \omega_0 \phi \; \rangle  =
\langle \;Ae , \omega_0\phi - \xi \; \rangle =
\langle \;\omega_2 Ae , \omega_0\phi - \xi \; \rangle\\
& =&
\langle \;\omega_2 A \omega_4 e , \omega_0\phi - \xi \; \rangle +
\langle \;\omega_2 A (1-\omega_4) e , \omega_0\phi - \xi \; \rangle
=:E+F. \end{eqnarray*}
Now from (\ref{supbsp}) and Lemma \ref{gla} we derive
$$ \vert E \vert \leq
\Vert \omega_2 A \omega_4 e \Vert_{-\alpha} 
\Vert \omega_0 \phi - \xi \Vert_{\alpha} \leq C\,  h 
\Vert e \Vert_{\alpha,\Gamma_*}
 \Vert \phi \Vert_{\alpha}
$$ 
$$ \vert F \vert \leq
\Vert \omega_2 A (1-\omega_4) e \Vert_{-\alpha}
\Vert \omega_0 \phi - \xi \Vert_{\alpha} \leq C\,  h
\Vert e \Vert_{\beta}
\Vert \phi \Vert_{\alpha}.
$$
Since $ (1-\omega_2) \omega_0=0 $ we have
$$
D= \langle \;(A\omega_0-\omega_0A) e, \omega_M \phi \;\rangle
$$
$$
=\langle \;(A\omega_0-\omega_0A) \omega_2 e, 
            \omega_M \phi \; \rangle +
\langle \;(A\omega_0-\omega_0A) (1-\omega_2) e, 
            \omega_M \phi \; \rangle 
$$
$$
=\langle \;(\tilde{A}\omega_0-\omega_0 \tilde{A}) \omega_2 e, 
            \omega_M \phi \; \rangle +
\langle \;-\omega_0A (1-\omega_2) e, 
            \omega_M \phi \; \rangle =: H_1+H_2
$$
Due to  Lemma \ref{com} and Lemma \ref{gla} we have the estimate
$$\vert H_1 \vert \leq 
\Vert (\tilde{A}\omega_0-\omega_0\tilde{A}) \omega_2 e
\Vert_{-\alpha} 
\Vert \omega_M \phi \Vert_{\alpha}
\leq C\,  \Vert e \Vert_{\alpha-1,\Gamma_*} \Vert\phi\Vert_{\alpha} $$
$$ \vert H_2 \vert \leq
\Vert \omega_0A (1-\omega_2) e \Vert_{-\alpha} 
\Vert\omega_M\phi\Vert_{\alpha} \leq C\,  
\Vert e \Vert_{\beta } \Vert\phi\Vert_{\alpha} $$
Putting together yields the desired estimate. \qed

\begin{lemma} \label{dec3}
 Under the assumptions (\ref{a1}-\ref{a3}) there holds
$$ \Vert G\tilde{\psi_h} - \tilde{\psi_h}\Vert_{\alpha}\leq C\,  
\left\{ h^{s-\alpha} \Vert \psi \Vert_{s,\Gamma_*} +
h \Vert  e \Vert_{\alpha,\Gamma_2} \right\}
$$
\end{lemma}
Proof: Due to the quasioptimal estimate for the error of the Galerkin solution
and to  Lemma \ref{sup} we have
\begin{equation} \label{hi0}
 \Vert G\tilde{\psi_h} - \tilde{\psi_h}\Vert_{\alpha} 
\leq C\,  \inf_{\xi \in S_h^r} \Vert \tilde{\psi _h} - \xi \Vert_{\alpha}
\leq C\,  h^{s+1-\alpha} \Vert \psi_h \Vert_{s,\Gamma_1}.
\end{equation}
For arbitrary $\phi \in S_h^r$ there holds
\begin{equation} \label{hi1}
\Vert \psi_h \Vert_{s,\Gamma_1} 
= \Vert \omega_1 \psi_h \Vert_{s,\Gamma_1}
\leq \Vert \omega_1 (\psi_h -\phi) \Vert_{s,\Gamma_1}
+\Vert \omega_2 \psi -\phi \Vert_{s,\Gamma_1}
+\Vert \omega_2 \psi  \Vert_{s,\Gamma_1}
\end{equation}
We choose $\phi$ due to Lemma \ref{app} b) such that there  holds for 
$\alpha \leq t\leq s$
\begin{equation} \label{hi2}
\Vert \omega_2 \psi -\phi \Vert_{t} \leq C\,  h^{s-t} \Vert \psi \Vert_{s,\Gamma_*}.
\end{equation} 
Applying  (\ref{locinv})  for the first term in (\ref{hi1}), and  using  the inequality (\ref{hi2}) we derive from (\ref{hi1})
\begin{eqnarray*}
\Vert \psi_h \Vert_{s,\Gamma_1} 
&\leq& C\,  \left\{ h^{-s+\alpha} \Vert \psi_h-\phi \Vert_{\alpha,\Gamma_2}
+ \Vert \psi \Vert_{s,\Gamma_*} \right\}
\\
&\leq & C\,  \left\{ 
h^{-s+\alpha} \Vert \psi_h-\psi \Vert_{\alpha,\Gamma_2}+
h^{-s+\alpha} \Vert \psi-\phi \Vert_{\alpha,\Gamma_2}
+ \Vert \psi \Vert_{s,\Gamma_*} \right\}
\\
& \leq&  C\,  \left\{ 
h^{-s+\alpha} \Vert \psi_h-\psi \Vert_{\alpha,\Gamma_2}+
h^{-s+\alpha} \Vert \omega_2\psi-\phi \Vert_{\alpha,\Gamma_2}
+ \Vert \psi \Vert_{s,\Gamma_*} \right\}\\
&\leq &C\,  \left\{ 
h^{-s+\alpha} \Vert \psi_h-\psi \Vert_{\alpha,\Gamma_2}
+ \Vert \psi \Vert_{s,\Gamma_*} \right\},
\end{eqnarray*}
where we have used (\ref{hi2})  with $t=\alpha$. 
Together with (\ref{hi0}), the last inequality yields 
$$ \Vert G\tilde{\psi_h} - \tilde{\psi_h}\Vert_{\alpha}\leq C\,  
\left\{ h^{s+1-\alpha} \Vert \psi \Vert_{s,\Gamma_*} +
h \Vert  e \Vert_{\alpha,\Gamma_2} \right\},
$$
and hence the lemma is proved.   \qed

Next we give in Lemma \ref{erinen} local error estimates in the energy norm, which is the 
$H^{\alpha}$ norm.

\begin{lemma} \label{erinen}
Assume that the assumptions  (\ref{a1}-\ref{a3}) hold. Let $\beta<\alpha $ arbitrary but fixed.
$$\Vert \psi -\psi_h \Vert_{\alpha,\Gamma_0} \leq C\,  \left\{ 
h^{s-\alpha} \Vert \psi \Vert_{s,\Gamma_*} 
        + \Vert \psi -\psi_h \Vert_\beta \right\}$$
\end{lemma}
Proof:
Combining Lemma \ref{dec1} - \ref{dec3} yields
$$\Vert \psi -\psi_h \Vert_{\alpha,\Gamma_0} \leq C\,  \left\{
h^{s-\alpha} \Vert \psi \Vert_{s,\Gamma_*} + 
h    \Vert \psi -\psi_h \Vert_{\alpha,\Gamma_1} 
+    \Vert \psi -\psi_h \Vert_{\beta} \right\}$$
and applying this estimate to $\Vert \psi -\psi_h \Vert_{\alpha,\Gamma_1} $
we obtain
$$\Vert \psi -\psi_h \Vert_{\alpha,\Gamma_0} \leq C\,  \left\{
h^{s-\alpha} \Vert \psi \Vert_{s,\Gamma_*} + 
h^2  \Vert \psi -\psi_h \Vert_{\alpha,\Gamma_2} 
+    \Vert \psi -\psi_h \Vert_{\beta} \right\} $$
We repeat this procedure and get the estimate 
$$\Vert \psi -\psi_h \Vert_{\alpha,\Gamma_0} \leq C\,  \left\{
h^{s-\alpha} \Vert \psi \Vert_{s,\Gamma_*} + 
h^M  \Vert \psi -\psi_h \Vert_{\alpha,\Gamma_M} 
+    \Vert \psi -\psi_h \Vert_{\beta}\right\}. $$
Further we estimate with  (\ref{locinv}) the second term by the first and the third term
$$ \Vert \psi -\psi_h \Vert_{\alpha,\Gamma_M} \leq 
 \Vert \psi  \Vert_{\alpha,\Gamma_M} 
+\Vert \psi_h \Vert_{\alpha,\Gamma_M}
\leq \Vert \psi  \Vert_{\alpha,\Gamma_M} 
+\Vert \omega_M\psi_h \Vert_{\alpha,\Gamma}$$
$$ \leq 
 \Vert \psi  \Vert_{s,\Gamma_*} 
+h^{\beta-\alpha}\Vert \psi_h \Vert_{\beta,\Gamma_*}
\leq \Vert \psi  \Vert_{s,\Gamma_*} 
+h^{\beta-\alpha}\Vert \psi -\psi_h \Vert_{\beta,\Gamma_*}
+h^{\beta-\alpha}\Vert \psi  \Vert_{\beta,\Gamma_*}
$$ $$
\leq C\, 
h^{\beta-\alpha}\Vert \psi  \Vert_{s,\Gamma_*}
+h^{\beta-\alpha}\Vert \psi -\psi_h \Vert_{\beta}.
$$
Now we choose $M$ such that $M \geq- \beta+s$,  
$$ h^M \Vert \psi - \psi_h \Vert_{\alpha,\Gamma_M} \leq 
 h^{\beta-\alpha+M} \Vert \psi \Vert_{s,\Gamma_*}
+h^{\beta-\alpha+M} \Vert \psi -\psi_h \Vert_{\beta}
\leq  h^{s-\alpha} \Vert \psi \Vert_{s,\Gamma_*}
+h^{s-\alpha} \Vert \psi -\psi_h \Vert_{\beta}$$ 
and hence the proof is complete. \qed 

Finally applying with Lemma \ref{erinen} and the Aubin-Nitsche trick 
shows that the local estimates can also be derived for Sobolev norms 
of negative orders (see \cite{thanh}). \qed
\vspace{0.3cm}

Common examples of these operators are Symm's integral equation with 
logarithmic kernel $(\alpha =-1/2)$ and the hypersingular integral equation
$(\alpha =1/2)$, defined respectively as
\begin{equation} \label{vtilde}
V \psi (x) := 
-\frac{1}{\pi} \int_{{\Gamma}} \; \log \vert x-y \vert \; \psi (y) \; ds_y
\end{equation}
\begin{equation} \label{dtilde}
W \psi (x) := 
-\frac{1}{\pi} \frac{\partial}{\partial n_x} \int_{{\Gamma}} \;
\frac{\partial}{\partial n_y} \log \vert x-y \vert \; \psi (y) \; ds_y,
\end{equation}
where $ \frac{\partial}{\partial n_x}$ denotes the derivate operator in direction of the 
outside normal.

In the case $A=V$ assumption (\ref{bb}) is satisfied if $ \Gamma $  has  capacity $ cap(\Gamma) < 1 $  (which can always be achieved after an appropriate scaling of $ \Gamma $. \cite{HW})
The following results are shown in \cite{cost1,cost2}:

\begin{theorem} \label{global}
There is a meshwidth $h_0>0$ such that for $ 0<h<h_0$ the Galerkin equations
(\ref{galeq}) with $A=V$ and $\alpha=-1/2$ are uniquely solvable in $S^0_h(\Gamma)$. Moreover, there holds
for $-s_0-1/2<t\leq s<s_0-1/2$:
\begin{equation} \label{galest}
\Vert \psi_h -\psi \Vert_t \leq C h^{s-t}\Vert \psi \Vert_s
\end{equation}
where $s_0$ is given by 
\begin{equation} \label{s0}
s_0 := \min \left\{ \frac{\pi}{\beta_j},  \frac{\pi}{2\pi -\beta_j}, 
              |j=1,\ldots,J \right\}
\end{equation}
and takes values between $1/2$ and  $1$.
\end{theorem}

Therefore, the order of convergence in the $L^2$ norm and in the energy norm (i.e. 
$H^{-1/2}$ norm) are ${\cal O}(h^{s_0-1/2-\epsilon}$) and 
{\cal O}$(h^{s_0-\epsilon}$), 
respectively, with $\epsilon >0$ arbitrary and $s_0$ from (\ref{s0}). The highest order 
 {\cal O}($h^{2(s_0-\epsilon)}$) is obtained in the $H^{-s_0-1/2+\epsilon}$ norm.

\begin{remark}  \label{rem1}
 If $\Gamma$ is the boundary of a L-shaped domain
with interior angles  $\frac{\pi}{2}$ and $\frac{3\pi}{2}$ we have for 
$s_0=2/3$, in (\ref{s0}). Hence we expect global convergence rates 
in the $L_2$ norm of order $O(h^{1/6-\epsilon})$ and in the $H^{-1/2}$ norm 
of order $O(h^{2/3-\epsilon})$. The highest rate 
$O(h^{4/3-2\epsilon})$ is obtained in the $H^{-7/6-\epsilon}$ norm.
\end{remark} 
\begin{remark}  \label{rem2}
In order to illustrate the results of Theorem \ref{theo} 
we consider the convergence rates
for the Galerkin method with piecewise constant trial functions
and the case $A=V$. 
Due to Theorem \ref{theo} we expect local 
convergence rates of order $O(h^\gamma)$ with 
$\gamma=\min\{1, 2(s_0-\epsilon)-1/2\}$ for the $L^2$ norm  and with
$\gamma=\min\{3/2, 2(s_0-\epsilon)\}$ for the $H^{-1/2}$ norm. 
More precisely let $\Gamma$ be an L-shaped polygon. In this case we expect
the local convergence in the $L^2$ norm to be of order $O(h^{5/6-\epsilon})$
and in the $H^{-1/2}$ norm of order 
$O(h^{4/3-\epsilon})$.
Thus we have higher local convergence than global convergence (see Remark \ref{rem1}).
\end{remark} 
\begin{remark}  \label{rem3}
For the hypersingular operator W a result corresponding to  Theorem \ref{global} can be found in 
\cite {cost3}.
\end{remark}

\section{Numerical Experiments} \label{nummm}

In this section we consider 
$V_\Gamma \psi =f$ on an L-shaped polygon $\Gamma$, 
whose sides are of length $1$ and $2$.  
The Galerkin method with piecewise constant trial functions (i.e. r=1)
is applied, as right hand side 
we choose $ f= (I+K)g $ for  $g=\Im (z^{2/3})$ where $K$ 
denotes the double layer potential operator. 
$$
K g (x) := -\frac{1}{\pi} \int_{{\Gamma}} \; 
\frac{\partial}{\partial n_y} \log \vert x-y \vert g (y) \; ds_y
$$ 
Here the exact solution  $\psi = \frac{\partial }{\partial_{n}} g$ 
contains singularities due to 
(\ref{s0}). In Table 1 the Galerkin error $\Vert \psi -\psi_N \Vert_{0,\Gamma_0}$
and experimental convergence rates are given with $h=N/8$.
In the headline of the table $a=0.02$ denotes that 
the Galerkin error is computed on $\Gamma_0$ which 
consists of straight subsegments which have 
a distance of  $0.02$ from the corners. 
Note for the global error (which corresponds to $a=0$) 
we achieve the expected order    
of convergence $O(h^{1/6-\epsilon})$, see (\ref{galest}). 
Theorem \ref{theo} predicts local convergence rates
of order $O(h^{5/6-\epsilon})$ due to the second term on the right hand side in (\ref{inn})
since $\beta= -7/6 +\epsilon$ and $\sigma = -1/2$ yields 
$h^\sigma \Vert \psi -\psi_h \Vert_{\beta,\Gamma} 
\leq C \, h^{5/6-\epsilon}$. The
numerical results given in Table~1 show we have a better convergence rate of $O(h)$
which seems to be due to the first term in the right hand side
in (\ref{inn}) (Note $t=0, s\leq 1=r$). 
This can be seen as follows.  
Due to Remark \ref{rem2} the local error
is controlled by the highest possible order of the global error 
in Sobolev norms of negative order. 
However, the order $O(h)$ is only obtained
if $8/N=h<a$, i.e. if the subsegment where the local error is evaluated
has a distance from the corner which is greater than the meshsize $h$.

\begin{center}

\vspace*{0.2cm}

\footnotesize
\begin{tabular}{@{}|r||c|c||c|c||c|c||c|c|}
\hline
  & \multicolumn{2}{|c||}{ $a=0$} &\multicolumn{2}{|c||}{ $a=0.02$} &
    \multicolumn{2}{|c||}{ $a=0.07$} &\multicolumn{2}{|c| }{ $a=0.15$} \\
\hline
N &  $\Vert \psi -\psi_N \Vert_{0,\Gamma_0}$ & EOC 
  &  $\Vert \psi -\psi_N \Vert_{0,\Gamma_0}$ & EOC
  &  $\Vert \psi -\psi_N \Vert_{0,\Gamma_0}$ & EOC 
  &  $\Vert \psi -\psi_N \Vert_{0,\Gamma_0}$ & EOC \\ \hline
 8 &  8.26 -1 & & 4.64 -1 & & 3.56 -1 & & 3.16 -1 & \\
  & & 0.16 & & 0.31 & & 0.18 & & 0.03 \\
 16&  7.37 -1 & & 3.74 -1 & & 3.13 -1 & & 3.09 -1 & \\
  & & 0.17 & & 0.34 & & 0.19 & & 0.30 \\
 32&  6.56 -1 & & 2.95 -1 & & 2.74 -1 & & 2.50 -1 & \\
  & & 0.17 & & 0.25 & & 0.26 & & 2.25 \\
 64&  5.84 -1 & & 2.47 -1 & & 2.29 -1 & & 5.27 -2 & \\
  & & 0.17 & & 0.18 & & 2.09 & & 1.26 \\
128&  5.21 -1 & & 2.19 -1 & & 5.40 -2 & & 2.20 -2 & \\
  & & 0.17 & & 0.34 & & 1.43 & & 1.00 \\
256&  4.64 -1 & & 1.73 -1 & & 2.00 -2 & & 1.10 -2 & \\
  & & 0.17 & & 2.39 & & 0.99 & & 1.01 \\
512&  4.13 -1 & & 3.29 -2 & & 1.01 -2 & & 5.48 -3 & \\
  & & 0.17 & & 1.06 & & 0.97 & & 0.97 \\
1024&  3.68 -1 & & 1.58 -2 & & 5.15 -3 & & 2.80 -3 & \\ \hline
\end{tabular}

\vspace{1cm}

\normalsize
Table 1: Local convergence rates and experimental orders of convergence 
\end{center}

\section{The K-operator} \label{kop}

In this section we consider a post processing method
which improves the order of local convergence. This better 
approximation is constructed by averaging the values of the Galerkin 
solution, using the K-operator. 
The action of the K-operator is given by the convolution of the 
Galerkin solution with a special spline.

The K-operator was proved to be an effective post processing
method in the finite element environment \cite{br,tho,whal}. 
In the boundary element literature the K-operator is applied 
to Galerkin solutions of strongly elliptic 
integral equations
on closed and open smooth curves by Tran
\cite{thanh,tran2,tran3} and on open surfaces  
by Stephan und Tran 
\cite{thste}.

In the following we apply the K-operator to the Galerkin solution with piecewise
constant trial functions for the integral equation $V_\Gamma \psi=f$,
where $\Gamma$ is a polygon, see ({\ref{vtilde}).
The results can be extended to other pseudo-differential operators, if 
global error estimates in negative norms are available.
To simplify the notation we define the straight line segment $\Gamma_j$ 
of the polygon $\Gamma $
 with the interval $I=[-1,1]$. Assume that the smooth closed curve 
$\tilde{\Gamma}$ containing $\Gamma_j$ is choosen in such a way, that 
the interval $I'=[-2,2]$ is contained in $\tilde{\Gamma}$, i.e. we have
$\Gamma_j = I \subset I' \subset \tilde{\Gamma}$.

For the application of the K-operator we assume that the mesh
is uniform on the subsegment 
$\Gamma_*$. 
Then there exists an $h_0>0$, such that for any $h\in (0,h_0]$, for
 $i=0,1,\ldots ,M-1$,
and for $j=1,2$:
$$ T^j_{\pm h}\phi \in S^{r,0}_h(\Gamma_{i+1}) \qquad \forall  
\phi \in S^{r,0}_h(\Gamma_i)$$
where $T_h$ denotes the translation operator $T_hv(x)=v(x+h)$ and 
$T_h^j$ for $j\geq 2$ is defined by
$T_h^jv=T_h(T_h^{j-1}v)$.

The K-operator acting on $\phi_h$ is defined by the convolution of $\phi_h$
with a function $K_h$ defined as a linear combination of B-splines.
For arbitrary but fixed $l,q\in \mathbb{N} $ we define 
for $0<h<1$ 
$$ K_h(x):=\frac{1}{h} K_q^l(\frac{x}{h})$$
with
$$K_q^l(x) = \sum_{j=-(q-1)}^{q-1}\; k_j \psi^{(l)}(x-j)$$
where $\psi^{(l)}$ is given by convoluting  $(l-1)$ times $\chi$, where
$$\chi(x) =\left\{ \begin{array}{cl}
1 & \mbox{ if } -1/2<x<1/2, \\
1/2 & \mbox{ if } x\pm 1/2,\\
0 & \mbox{ otherwise } 
\end{array} \right.
$$
Hence $\psi^{(l)}$ is the  B-spline of order $l$ and symmetric 
about 0, with support 
$[-1/2,1/2]$.

Here  $k_j$ for $ j=-(q-1),\ldots ,q-1$ is chosen such that
\begin{equation} \label{defkj}
 \int_{-\infty}^\infty\; K_q^l(\frac{x}{h})x^i \;dx=
\left\{
\begin{array}{cll}
 1 &  \mbox{f\"ur} & i=0\\
 0 &  \mbox{f\"ur} & i=1,\ldots ,2q-1
\end{array} \right.
\end{equation}
and satisfies the symmetry condition $k_j=k_{-j}$ for $j=1,\ldots ,q-1$.
The K-operator is defined in such a way, that it reproduces polynomials 
(up to degree $2q-1$) under convolution.
We recall some properties of the K-operator from \cite{br,tran3} and refer
for further details to these papers.

\begin{lemma} \label{lemk} \cite{tran3}
For any $ i=0,\ldots , M-1$, there exists $c>0$ and $h_0>0$, such that
for any $v \in H^{s}(\Gamma_{i+1})$ with $0\leq s \leq 2q$:
$$ \Vert K_h *v -v \Vert_{0,\Gamma_i} \leq c h^{s} \Vert v \Vert_{s,\Gamma_{i+1}}$$
\end{lemma}

Another interesting property of $K_h$ is that the differential operator
when applied to 
$K_h* v$ is changed to the central differential operator applied
to a somewhat
similar function. More precisely, letting
\begin{equation} \label{defd}
\partial_h v(x) = \frac{1}{h}\{ v(x+h)-v(x)\}=\frac{1}{h} \{T_h v(x) -v(x) \}
\end{equation}

we have the following
\begin{lemma} \label{lemp} \cite{tran3}
For any $j=0,1,\ldots ,l$ and $i=0,1,\ldots ,M-1$ we have 
$$ \Vert D^j(K_h * v)\Vert_{s,\Gamma_i}
\leq c \Vert \partial^j_h v\Vert_{s,\Gamma_{i+1}}$$
\end{lemma}

\begin{lemma} \label{lemn} \cite{br}
Let $\tau$ be a non-negative integer. Then for any 
$i=0,1,\ldots ,M-1$ there exists a constant $C$, such that:
$$ \Vert v \Vert_{0,\Gamma_i}
\leq C\,  \sum_{j=0}^\tau \Vert D^j v \Vert_{-\tau,\Gamma_{i+1}}
$$
\end{lemma}

To apply the K-operator to the Galerkin solution 
given on a polygon we define for any function $v$
defined on $\Gamma$ the function $\tilde{v}$ to be the restriction
of $v$ on $\Gamma^j$ extended to $\mathbb{R} $ by zero.
Further we introduce the cut-off function
$$\omega_* \in C_0^\infty(\Gamma_*) \quad \omega_* \equiv 1 \; \mbox{on}\;
\Gamma_M
$$
and define the acting of the K-operators on $v$ by
\begin{equation} \label{actk}
 K_h(v):=K_h* (\omega_* \tilde{v}). 
\end{equation}

\begin{theorem} \label{koper}
Let $l=2$ and $q=2$. For the solution $\psi$ of $V\psi=f$ there holds
$\psi \in H^{3}(\Gamma_*)$. Further let $\psi_h$ be the Galerkin solution
given by (\ref{galeq}) with $r=1$ and let $K_h(\psi_h)$ defined as in 
(\ref{actk}). Then there exists a constant 
$C$ independent of $h$
\begin{equation} \label{kerr} \Vert K_h( \psi_h) - \psi \Vert_{0,\Gamma_0} \leq C\,  \{h^{3}
\Vert \psi \Vert_{3,\Gamma_*} + h^{2(s_0-\epsilon)}
\Vert \psi  \Vert_{s_0-\epsilon,\Gamma}\}
\end{equation}
where $s_0$ is defined by (\ref{s0}).
\end{theorem}

Proof:
Using the triangle inequality and $\omega_* \psi = \tilde{\psi} $  yields
$$\Vert K_h( \psi_h) - \psi \Vert_{0,\Gamma_0}\leq
\Vert K_h*(\omega_* \tilde{\psi_h} - \omega_* \tilde{\psi}) \Vert_{0,\Gamma_0} +
\Vert K_h*(\omega_* \psi) - \omega_* \tilde{\psi} \Vert_{0,\Gamma_0}.$$
We estimate both terms seperately. 
Due to Lemma \ref{lemk} we have for the second term 
(Note: $q=2$)
\begin{equation} \label{sec}
 \Vert K_h*(\omega_* \psi) - \omega_* \tilde{\psi} \Vert_{0,\Gamma_0}
\leq C\,  h^{3} \Vert \psi \Vert_{3,\Gamma_1}.
\end{equation}
For the first term we derive with Lemma \ref{lemn} and \ref{lemp}  
\begin{eqnarray*}
\Vert K_h*(\omega_* \tilde{\psi_h} - \omega_*\tilde{\psi}) \Vert_{0,\Gamma_0}
&\leq & C\,  \sum_{j=0}^{2} \;
\Vert D^j K_h* (\omega_*\tilde{\psi_h} -\omega_*\tilde{\psi}) \Vert_{-2,\Gamma_1}\\
&\leq & C\,  \sum_{j=0}^{2} \;
\Vert \partial_h^j  (\omega_*\tilde{\psi_h} -\omega_*\tilde{\psi}) \Vert_{-2,\Gamma_2}.
\end{eqnarray*}
and since $T_{\pm 2h}(\Gamma_2) \subset \Gamma_*$
$$
\Vert \partial_h^j  (\omega_*\tilde{\psi_h} -\omega_*\tilde{\psi}) \Vert_{-2,\Gamma_2} =
\Vert \partial_h^j  (\tilde{\psi_h} -\tilde{\psi}) \Vert_{-2,\Gamma_2}.
$$
First we proof the auxiliary identity
\begin{equation} \label{hilf}
 \langle \; \tilde{V} \partial_h^j (\tilde{\psi}_h-\tilde{\psi}),\phi \; \rangle =0
 \qquad
\forall \phi \in S^{1,0}_h(\Gamma_2),
\end{equation}
where $\tilde{V}$ denotes the extension of $V$ to $\tilde{\Gamma}$ and
$\langle \cdot , \cdot \rangle$ denotes the duality in $L^2(\tilde{\Gamma})$.
Now (\ref{hilf}) can be derived from the definition of $\partial_h$ in 
(\ref{defd}) and the  Galerkin equations (\ref{a3}). 
$$
\langle \; \tilde{V} \partial_h^j (\tilde{\psi}_h-\tilde{\psi}),\phi \; \rangle
= \frac{1}{h^j} \sum_{i=0}^j \; { j \choose i}
\langle \; \tilde{V} T_h^i (\tilde{\psi}_h-\tilde{\psi}),\phi \; \rangle $$ 
$$
= \frac{1}{h^j} \sum_{i=0}^j \; {j \choose i}
\langle \;  T_h^i \tilde{V} (\tilde{\psi}_h-\tilde{\psi}),\phi \; \rangle
= \frac{1}{h^j} \sum_{i=0}^j \; {j \choose i}
\langle \;  \tilde{V} (\tilde{\psi}_h-\tilde{\psi}), T_{-h}^{i}\phi \; \rangle
=0.
$$
For exchanging $ \tilde{V} T_h^i$ by  $ T_h^i \tilde{V}$ we need
that the interval $[-2,2]$ is contained in $\tilde{\Gamma}$.
Note that for any $\phi \in S^{1,0}_h(\Gamma_i)$ we have $T^j_{\pm h}
\phi \in S_h^{1,0}(\Gamma_{i+1})$.
From  (\ref{hilf}) we see that  $\partial_h^j \psi_h$ is the ''interior''
Galerkin solution to $\partial_h^j \psi$. Now we apply the local
error estimates (\ref{inn}) with $\partial_h^j \psi-\xi $ instead of 
 $\psi$ and $\partial_h^j \psi_h - \xi$ instead of $\psi_h$, where  
$\xi \in S_h^{1,0}(\Gamma_{j+5})$ is arbitrary.
$$\Vert \partial_h^j (\psi_h-\psi) \Vert_{-2, \Gamma_2} \leq
C\,  \{h^{5/2-\epsilon }
\Vert \partial_h^j \psi -\xi \Vert_{1/2-\epsilon,\Gamma_3}
+\Vert \partial_h^j (\psi_h-\psi) \Vert_{-\beta,\Gamma} \}$$
\begin{equation} \label{xi}
\leq C\,  \{ h^{5/2-\epsilon}
\Vert \partial_h^j \psi - \xi \Vert_{1/2-\epsilon,\Gamma_3}
+ \Vert \psi_h -\psi \Vert_{\beta+j}\}.
\end{equation}
Further we obtain
$$ \Vert \partial_h^j \psi - \xi \Vert_{1/2-\epsilon,\Gamma_3} =
 \Vert \omega_3 \partial_h^j \psi - \xi \Vert_{1/2-\epsilon,\Gamma_3}
\leq  \Vert \omega_3 \partial_h^j \psi - \xi \Vert_{1/2-\epsilon,\Gamma}$$
Due to Lemma \ref{app} we choose $\xi$ such that 
 $$\Vert \omega_3 \partial_h^j \psi - \xi \Vert_{1/2-\epsilon,\Gamma}
  \leq C\,  h^{1/2+ \epsilon}\Vert \partial^j_h \psi \Vert_{1,\Gamma_4}
 \leq C\,  h^{1/2+ \epsilon}\Vert \psi \Vert_{j+1,\Gamma_4}$$
\begin{equation} \label{par}
 \leq C\,  h^{1/2+ \epsilon}\Vert \psi \Vert_{3,\Gamma_4}.
\end{equation}
where we used $j\leq 2$.
Now putting together  (\ref{par}), (\ref{xi}) and (\ref{sec}) yields
\begin{equation} \label{last} 
\Vert K_h( \psi_h) - \psi \Vert_{0,\Gamma_0}
\leq C\,  \{ h^{3} \Vert \psi \Vert_{3,\Gamma_*} +
\Vert \psi_h -\psi  \Vert_{\beta +2,\Gamma}\}, \end{equation}
and we obtain the desired estimation if we choose
$\beta$ such that 
$$ \Vert \psi_h -\psi  \Vert_{\beta +2,\Gamma} \leq
\Vert \psi_h -\psi  \Vert_{-s_0+\epsilon,\Gamma}
\leq c h^{2(s_0-\epsilon)} \Vert \psi \Vert_{s_0-\epsilon,\Gamma},
$$
where we have used Theorem \ref{global} for the last inequality. \qed

\begin{remark}  \label{rem4}
Again we will look more closely at the orders of convergence
if we consider the equation $V \psi =f$ on an L-shaped boundary. Let  
Galerkin's method with piecewise constant trial functions be computed for this 
equation. Then we expect due to Theorem \ref{koper} a local convergence of order
$O(h^{4/3-\epsilon})$, whereas due to Theorem \ref{theo} we obtain 
$O(h^{5/6-\epsilon})$ without application of the K-operator, see Remark \ref{rem2}. 
\end{remark}
Now as conclusion we describe a strategy how to improve the local convergence.
The local order of covergence is controlled  by 
the highest possible order for the global error in negative norms, i.e.
by the second term on the right hand side of (\ref{kerr}). This error
can be estimated using Theorem \ref{global} with $t$ negative.
In order to achieve the order of $O(h^3)$, which is predicted
by the first term of (\ref{kerr}), a partially graded mesh has to be used
in order to improve the global error in (\ref{last}) which yields the second term in 
(\ref{kerr}).
On the other hand for the application of the K-operator we need a uniform
mesh only on the subsegment $\Gamma_*$. This allows a graded mesh
near the corners. Using a
uniform mesh on $\Gamma_*$ and a graded mesh near the corners, higher
orders of global converence can be achieved with the results of \cite{peter}. 
Applying the K-operator to the Galerkin solution on such a combined
mesh yields  local convergence of order
$O(h^3)$ (cf. \cite{tran3} for the case of an interval $\Gamma$).

\end{document}